\documentclass[11pt]{article}
\usepackage[utf8]{inputenc}
\usepackage[T1]{fontenc}
\usepackage[a4paper,bindingoffset=0cm,head=15pt,body={17cm,25cm},top=2.5cm]{geometry}
\usepackage{amsmath,amssymb}
\usepackage{graphicx}
\usepackage{color}
\usepackage{cite}
\usepackage{ifthen}
\usepackage{titlesec}
\usepackage{fancyhdr}
\usepackage{booktabs}
\usepackage[english]{babel}
\usepackage{hyperref}
\usepackage{titlecaps}
\usepackage{xparse}

\titleformat{\section}{\normalfont\bfseries}{\thesection}{2ex}{}
\titlespacing*{\section}{0pt}{2.5ex}{1ex}
\titleformat{\subsection}{\normalfont\bfseries}{\thesubsection}{2ex}{}
\titlespacing*{\subsection}{0pt}{2.5ex}{1ex}
\newtheorem{definition}{Definition}

\newtheorem{theorem}[definition]{Theorem}

\newcounter{authorcounter}
\newcounter{institutecounter}

\newcommand{\wavesauthorlist}{}
\newcommand{\wavesaddresslist}{}
\newcommand{\wavesemail}{}
\newcommand{\wavesfootnotes}{}
\newcommand{\wavesauthorpre}{}

\def\theNumberTest#1{%
  \if\relax\detokenize\expandafter{\romannumeral-0#1}\relax
    true%
  \else
    false%
  \fi
}

\NewDocumentCommand{\wavesspeaker}{ O{} O{} m m}{%
    \ifthenelse{\value{authorcounter} > 1}{%
      \renewcommand{\wavesauthorpre}{, }%
    }{%
      \renewcommand{\wavesauthorpre}{}%
    }%
    \ifx\relax#1\relax
      \renewcommand{\wavesfootnotes}{}
    \else   
      \renewcommand{\wavesemail}{$^\ast$Email: #1}%
      \renewcommand{\wavesfootnotes}{, \ast}
    \fi
    \ifthenelse{\equal{\theNumberTest{#4}}{true}}{%
      \edef\wavesauthorlist{\wavesauthorlist%
        \wavesauthorpre{}\underline{#3}$^{#2%
        }$%
      }%
    }{%    
      \edef\wavesauthorlist{\wavesauthorlist%
        \wavesauthorpre\underline{#3}$^{\arabic{authorcounter}%
        }$%
      }
      \edef\wavesaddresslist{\wavesaddresslist% 
        \par%
        $^{\arabic{authorcounter}}$#4%
      }%      
      \stepcounter{authorcounter}%
      \stepcounter{institutecounter}
    }%

    \ifx\relax#2\relax
          \edef\wavesauthorlist{\wavesauthorlist%
        \wavesauthorpre$^{\wavesfootnotes%
        }$%
      }
    \else
    \ifthenelse{\equal{\theNumberTest{#2}}{true}}{%
      \edef\wavesauthorlist{\wavesauthorlist%
        \wavesauthorpre{}$^{,#2\wavesfootnotes%
        }$%
      }%
    }{%    
      \edef\wavesauthorlist{\wavesauthorlist%
        \wavesauthorpre$^{,\arabic{institutecounter}\wavesfootnotes%
        }$%
      }
      \edef\wavesaddresslist{\wavesaddresslist% 
        \par%
        $^{\arabic{institutecounter}}$#2%
      }%      
      \stepcounter{institutecounter}%
    }%
    \fi
  \ignorespaces
}

\NewDocumentCommand{\wavesauthor}{ O{} O{} m m}{%
    \ifthenelse{\value{authorcounter} > 1}{%
      \renewcommand{\wavesauthorpre}{, }%
    }{%
      \renewcommand{\wavesauthorpre}{}%
    }%
    \ifx\relax#1\relax
      \renewcommand{\wavesfootnotes}{}
    \else   
      \renewcommand{\wavesemail}{$^\ast$Email: #1}%
      \renewcommand{\wavesfootnotes}{, \ast}
    \fi%
    \ifthenelse{\equal{\theNumberTest{#4}}{true}}{%
      \edef\wavesauthorlist{\wavesauthorlist%
        \wavesauthorpre{}#3$^{#4%
        }$%
      }%
    }{% 
      \edef\wavesauthorlist{\wavesauthorlist%
         \wavesauthorpre{}#3$^{\arabic{institutecounter}%
         }$%
      }
      \edef\wavesaddresslist{\wavesaddresslist%
        \par%
        $^{\arabic{institutecounter}}$#4%
      }%
      \stepcounter{authorcounter}
      \stepcounter{institutecounter}%
    }%

    \ifx\relax#2\relax
          \edef\wavesauthorlist{\wavesauthorlist%
        $^{\wavesfootnotes%
        }$%
      }
    \else
    \ifthenelse{\equal{\theNumberTest{#2}}{true}}{%
      \edef\wavesauthorlist{\wavesauthorlist%
        $^{,#4\wavesfootnotes%
        }$%
      }%
    }{%    
      \edef\wavesauthorlist{\wavesauthorlist%
        $^{,\arabic{institutecounter}\wavesfootnotes%
        }$%
      }
      \edef\wavesaddresslist{\wavesaddresslist% 
        \par%
        $^{\arabic{institutecounter}}$#2%
      }%      
      \stepcounter{institutecounter}%
    }%
    \fi
  \ignorespaces
}

\fancypagestyle{plain}{%
  \fancyhead{}
}

\newenvironment{wavespaper}[3]{%
  \renewcommand{\wavesauthorlist}{}%
  \renewcommand{\wavesemail}{}%
  \setcounter{authorcounter}{1}%
  \setcounter{institutecounter}{1}%
     #2
  \twocolumn[
    \begin{center}
     \bfseries
     #1
     \bigskip

     \wavesauthorlist
     \mdseries
     \smallskip

     \wavesaddresslist
     \smallskip
 
     \wavesemail
    \end{center}%
  ]
}{%
}

\newcommand{\EW}{{\mathrm{EW}}}
\newcommand{\IC}{{\mathbb C}}
\newcommand{\IN}{{\mathbb N}}
\newcommand{\IR}{{\mathbb R}}
\newcommand{\IS}{{\mathbb S}}
\newcommand{\IT}{{\mathbb T}}
\newcommand{\bb}{{\mathbf b}}
\newcommand{\bd}{{\mathbf d}}
\newcommand{\be}{{\mathbf e}}
\newcommand{\bn}{{\mathbf n}}
\newcommand{\bp}{{\mathbf p}}
\newcommand{\bu}{{\mathbf u}}
\newcommand{\bv}{{\mathbf v}}

\newcommand{\bx}{{\mathbf x}}
\newcommand{\by}{{\mathbf y}}

\newcommand{\mC}{{\mathbf C}}
\newcommand{\mD}{{\mathbf D}}
\newcommand{\mI}{{\mathbf I}}
\newcommand{\mU}{{\mathbf U}}
\newcommand{\mV}{{\mathbf V}}
\newcommand{\mSigma}{{\mathbf\Sigma}}
\newcommand{\dn}{{\partial_\bn}}
\newcommand{\dnK}{{\partial_{\bn_K}}}
\newcommand{\calF}{{\mathcal F}}
\newcommand{\calT}{{\mathcal T}}
\newcommand{\deO}{{\partial\Omega}}
\newcommand{\deK}{{\partial K}}
\newcommand{\test}{^{\mathrm{test}}}
\newcommand{\trial}{^{\mathrm{trial}}}

\begin{document}

\begin{wavespaper}{
Trefftz methods with evanescent plane waves
}{%
% Put the list of authors here, using \wavesauthor or \wavesspeaker commands
  \wavesspeaker[andrea.moiola@unipv.it]{Andrea Moiola}{Department of Mathematics ``F.~Casorati'', University of Pavia, Pavia, Italy}
  \wavesauthor{Nicola Galante}{Sorbonne Université, Université Paris Cité, CNRS, INRIA, Laboratoire Jacques-Louis Lions, LJLL, EPC ALPINES, Paris, F-75005, France}
\wavesauthor[][]{Emile Parolin}{2}
}

\section*{Abstract}
Classical Trefftz methods approximate Helm\-holtz solutions using propagative plane waves and are subject to strong numerical instabilities.
Evanescent plane wave bases can substantially mitigate this phenomenon.
We propose a simple recipe to select such basis functions.
We show that the numerical results obtained by the Ultraweak Variational Formulation (UWVF) greatly improve thanks to this choice.
More details and examples will soon be available in \cite{GaMoPa2026}.

\smallskip

\noindent\textbf{Keywords:} Evanescent plane wave, Trefftz discontinuous Galerkin, Ultraweak variational formulation, Oversampling, Regularisation.

\section{Trefftz methods and plane waves}

Exactly 100 years ago, at the International Con\-gress of Applied Mechanics that took place in Z\"urich in 1926, Erich Trefftz proposed the approximate solution of a Laplace bound\-ary value problem (BVP) by piecewise-harmonic func\-tions \cite{Trefftz1926}.
Decades later, his work gave the name to a wide class of numerical methods applicable to any linear PDE: \textbf{Trefftz methods} approximate BVP solutions with discrete spaces made of piecewise solutions of the underlying PDE.
Trefftz methods often take the form of discontinuous Galerkin schemes (DG).

Trefftz methods have been extensively used to approximate the \textbf{Helm\-holtz equation}~\cite{TrefftzSurvey}
\begin{equation}\label{eq:Helmholtz}
\Delta u+\kappa^2u=0
\end{equation}
on a subdomain of $\IR^n$, where $\kappa>0$ is the wavenumber.
Classical finite element and DG methods approximate solutions of \eqref{eq:Helmholtz} by piecewise polynomials.
Since there is no polynomial solution of \eqref{eq:Helmholtz} other than the constant zero, Trefftz methods for this PDE resort to non-pol\-y\-no\-mi\-al bases.
Most Trefftz schemes (e.g.\ \cite{CED98,HiMoPe11,HMK02}) for \eqref{eq:Helmholtz} use \textbf{propagative plane wave} (PPW) basis functions:
\begin{equation}\label{eq:PPW} 
\bx\mapsto e^{\imath\kappa\bd\cdot\bx} \quad\text{with}\quad \bd\in\IS^{n-1},
\end{equation}
where $\IS^{n-1}:=\{\bv\in\IR^n: |\bv|=1\}$.
Discrete spaces spanned by PPWs with different propagation directions $\bd$ approximate general Helm\-holtz solutions with fast convergence rates \cite{MoHiPe11}.
However, it is widely acknowledged that PPW-based Trefftz implementations suffer from strong numerical instabilities.

\section{A negative PPW approximation result}
The instabilities of PPW-based Trefftz schemes have been described and addressed in several ways, mostly in terms of near-dependence of the PPWs, ill-conditioning of the linear systems, and linear algebra techniques (e.g.\ \cite[\S4.2]{HMK02}, \cite{BaBeDiTo2021}, \cite{CoNi2025}, \cite[\S4.3]{TrefftzSurvey}).
We follow a different approach, inspired by \cite{AdHu2019}, and report the result of \cite[Lem\-ma 4.2]{Parolin2023} (2D) and \cite[Lemma~4.4]{Galante2025} (3D).

\begin{theorem}\label{thm:PPW}
Let $B$ be the unit ball in $\IR^n$, $n\in\{2,3\}$, and $\kappa>0$.
Take any 
\begin{itemize}
\item target relative accuracy $0<\delta<1$,
\item large threshold value $M>1$.
\end{itemize}
Then, there exists $u\in C^\infty(\IR^n)$, solution of \eqref{eq:Helmholtz} on $\IR^n$, normalised as $\|u\|_{H^1(B)}=1$, such that, if $u$ is approximable by a linear combination of a finite number of PPWs with accuracy $\delta$, i.e.
$$\bigg\|u-\sum_{p=1}^P\mu_p e^{\imath\kappa\bd_p\cdot\bx}\bigg\|_{H^1(B)} 
\le \delta
\quad \substack{\displaystyle\text{for some }P\in\IN,\\
\displaystyle\boldsymbol\mu\in\IC^P,\\ 
\displaystyle\{\bd_p\}_{p=1}^P\subset\IS^{n-1},}
$$
then the Euclidean norm of the coefficient vector is larger than $M$:
$\|\boldsymbol\mu\|_{\IC^P}\ge M$.
\end{theorem}
The function $u$ is simply a high-index circular or spherical wave.
The $H^1(B)$ norm in the statement could be replaced by any finite-order Sobolev norm.

Theorem~\ref{thm:PPW} shows that, for some target Helm\-holtz solutions, \textbf{accurate PPW approximations necessarily require huge coefficients}.
If the function we seek to approximate is the $u$ provided by the theorem for large $M$ (e.g.\ $M$ larger than the reciprocal of machine precision), then numerical cancellation is inevitable and there is no way to stably represent its PPW approximation in computer arithmetic.

This is the fundamental issue underlying all PPW-Trefftz instabilities.
Linear algebra techniques (such as preconditioning or matrix decompositions) can not help finding an approximation if an approximation that is representable in computer arithmetic does not exist.
Recall also that \cite{AdHu2019} shows that, roughly speaking, stable numerical approximations are possible if and only if small-coefficient approximation exists.

Theorem~\ref{thm:PPW} states that accurate approximations and small-coefficient representations are mutually exclusive for \textit{some} Helm\-holtz solutions (the ``evanescent modes'').
Other Helmholtz solutions can be stably approximated by PPWs, as witnessed by the numerous successful uses of Trefftz schemes \cite{TrefftzSurvey}.

\section{Evanescent plane waves}\label{s:EPW}
To deal with the PPW approximation instability, following \cite{Parolin2023,Galante2025} we propose the use of evanescent plane waves (EPW):
\begin{equation}\label{eq:EPW}
\bx\mapsto e^{\imath\kappa\bd\cdot\bx}, \quad \bd\in\IC^n,\quad \bd\cdot\bd=\sum_{j=1}^nd_j^2=1.
\end{equation}
EPWs allow complex propagation vectors $\bd$ and include PPWs as special case ($\bd\in\IR^n\subset\IC^n$).
The condition $\bd\cdot\bd=1$ is equivalent to 
\begin{equation}\label{eq:EPWcondition}
|\Re(\bd)|^2 - |\Im(\bd)|^2 = 1,\quad
\Re(\bd) \cdot \Im(\bd) = 0.
\end{equation}
An EPW is determined by three parameters:
\begin{align*}
\text{propagation direction }\ &\bp := \frac{\Re(\bd)}{|\Re(\bd)|} \in\IS^{n-1},\\
\text{evanescence direction }\ &\be := \frac{\Im(\bd)}{|\Im(\bd)|} \in\IS^{n-1},\\
\text{evanescence strength }\ &\eta := {|\Im(\bd)|} \in [0,+\infty).
\end{align*}
Note that $\be$ has to be orthogonal to $\bp$ by \eqref{eq:EPWcondition}, and that if $\eta=0$ then the EPW is a PPW and $\be$ is undefined.
We also introduce a second measure of the evanescence strength:
\begin{equation*}
    \zeta := {|\Re(\bd)|} = \sqrt{1+\eta^{2}} \in [1,+\infty).
\end{equation*}
Then the evanescent plane wave is
\begin{equation}\label{eq:EWy}
\EW_{[\bp,\be,\eta]}(\bx) := e^{\imath\kappa\bd\cdot\bx}=
e^{\imath \kappa \zeta \bp \cdot \bx}
e^{      -\kappa \eta  \be \cdot \bx}.
\end{equation}
This EPW oscillates and propagates in direction $\bp$ with apparent wavenumber $\kappa\zeta$, and decays with exponential rate $\kappa\eta$ in direction $\be$; see Figure~\ref{fig:Waves}.

\begin{figure}[htb]\centering
\includegraphics[width=0.3\linewidth,clip,trim=85 30 70 20,clip]{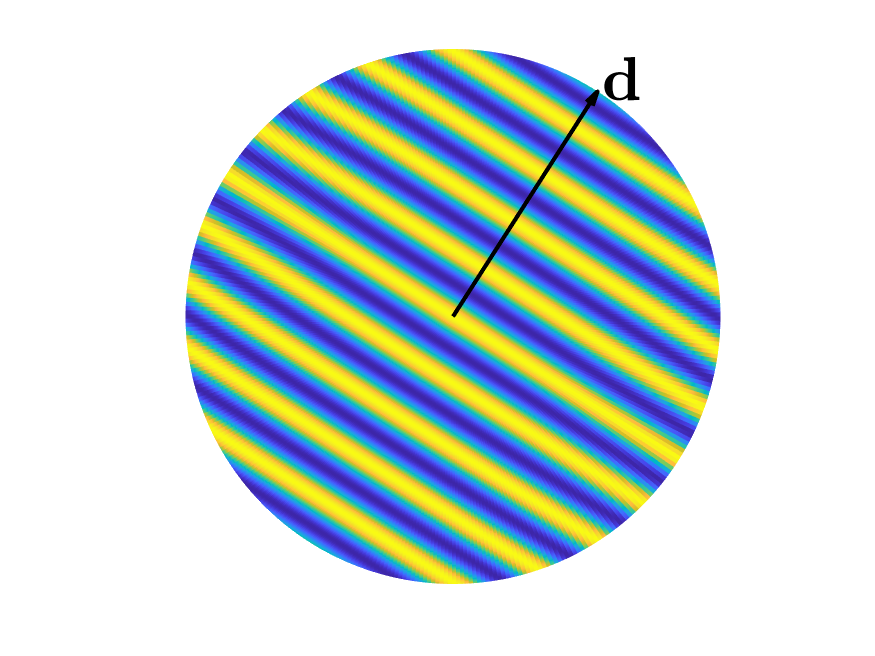}
\includegraphics[width=0.3\linewidth,clip,trim=0 0 0 0,clip]{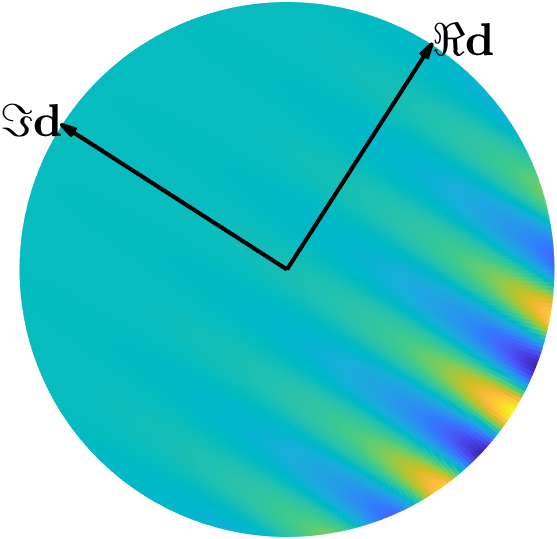}
\includegraphics[width=0.37\linewidth,clip,trim=170 180 35 130,clip]{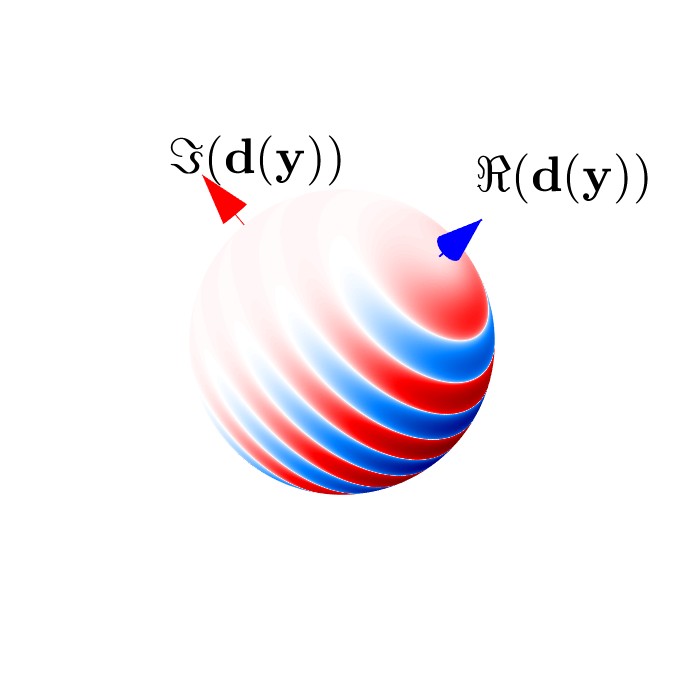}
\caption{A PPW and an EPW on the disc in $\IR^2$, an EPW on the ball in $\IR^3$. \label{fig:Waves}}
\end{figure}

Why do we expect linear combinations of the EPWs \eqref{eq:EPW} to give accurate small-coefficient approximations?
The reason is in \cite[Thm.~6.7]{Parolin2023} (2D) and \cite[Thm.~3.9]{Galante2025} (3D):
any Helmholtz solution $u\in H^1(B)$, with $B$ the unit disc/ball, is a continuous EPW superposition.
The coefficient of this superposition, measured in a weighted $L^2$ norm, is bounded by the $H^1(B)$ norm of $u$.
No such a representation can hold for PPWs.
This ``Herglotz transform'' is a step towards a convergence and stability theory for EPW-based schemes; further work is needed to extend it to more general domains and discrete EPW spaces.

\section{EPW selection}
Implementing a PPW-based Trefftz method, one typically selects a $P$-dimensional discrete space by choosing PPW directions $\{\bd_p\}_{p=1}^P$ that are roughly equispaced on the unit sphere $\IS^{n-1}$.

When implementing an EPW-based method, how should one choose the EPW parameters $\bp,\be,\eta$?
The propagation directions $\bp$ can be taken uniformly on $\IS^{n-1}$ and the evanescence directions $\be$ uniformly in $\{\be\in\IS^{n-1}:\be\cdot\bp=0\}\sim\IS^{n-2}$.
The choice of the evanescence strengths $\eta$ of the basis functions is less obvious but is key to ensure the quality of the method.

We describe in detail a simple concrete strategy to choose a discrete EPW space in 2D ($n=2$).
We refer to \cite{GaMoPa2026} for the 3D case.

In 2D, the propagation direction $\bp$ is pa\-ram\-e\-trised by the polar angle $\theta\in[0,2\pi)$ and the decay direction $\be$ by $\varphi\in\{\pm1\}$, as $\be$ can only be one of the two unit vectors orthogonal to $\bp$.
So \eqref{eq:EWy} reads
$$
\EW_{[\bp,\be,\eta]}(\bx) = 
e^{\imath\kappa\zeta\binom{\cos\theta}{\sin\theta}\cdot\bx} 
e^{-\kappa\eta\varphi\binom{-\sin\theta}{\cos\theta}\cdot\bx}.
$$

\subsection{EPW sampling recipe in 2D}\label{s:Recipe}
Assume we are given a ``DOF budget'' $P\in\IN$ and we want to construct a $P$-dimensional EPW Trefftz space on a polygonal domain $K\subset\IR^2$ (which in \S\ref{s:UWVF} will represent a mesh element).
We define a ``Fourier truncation level'' 
$$L:=\frac P4,$$
borrowing the name from \cite[\S7.3]{Parolin2023}, \cite[p.~462]{Galante2025}.
We sample $P$ points $\{\by_p=(\theta_p,\varphi_p,\xi_p)\}_{p=1}^P$ from the uniform probability measure on 
$$
Y:=[0,2\pi)\times\{\pm1\}\times[0,1].
$$
This probability is the Lebesgue measure scaled by $\frac1{2\pi}$ in the first component $\theta_p$, equal probability $\frac12$ for $\varphi_p=1$ and $\varphi_p=-1$, and the Lebesgue measure in the last component $\xi_p$.
We define the sampled evanescence parameters 
$$
\zeta_p:=\max\bigg\{1,\frac{2L}{\kappa\,\mathrm{diam}(K)}\xi_p\bigg\} 
\;\;\text{and}\;\;
\eta_p^2:=\zeta_p^2-1.
$$
The discrete Trefftz space is 
\begin{align}\nonumber
\IT_P(K):=&\mathrm{span}\bigg\{\frac{\EW_{\by_p}}{\|\EW_{\by_p}\|_{L^{\infty}(K)}}\bigg\}_{p=1}^{P},\; \text{with}
\\
\EW_{\by_p}(\bx):=&
e^{\imath\kappa\zeta_p\binom{\cos\theta_p}{\sin\theta_p}\cdot\bx} 
e^{-\kappa\eta_p\varphi_p\binom{-\sin\theta_p}{\cos\theta_p}\cdot\bx}.
\label{eq:EPWspace}
\end{align}
The norm \(\|\EW_{\by_p}\|_{L^{\infty}(K)}\) can be computed by evaluating the EPW at each vertex of the cell~\(K\).

\subsection{Comments on the sampling recipe}
The choices made in \S\ref{s:Recipe} might seem arbitrary: they are motivated in \cite{GaMoPa2026} from an asymptotic approximation of the density appearing in the Herglotz integral representation \cite{Parolin2023} (see also \cite[eq.~(5.8)]{Galante2025}) mentioned in \S\ref{s:EPW}.

The sampling of the set $Y$ can be done in several different ways: random, quasi-random (e.g.\ with Sobol sequences, as in~\S\ref{s:Numerics}), deterministic (e.g.\ equispaced in $\theta_p$ and/or $\xi_p$), with or without a tensor-product structure.

If $2L\le \kappa\,\mathrm{diam}(K)$, then all $\zeta_p=1$ and all basis functions are PPWs.
In this case the DOF budget $P$ is low, and all the approximation effort is devoted to the approximation of the propagative components of the solution.
For larger values of $P$, an increasing fraction $\frac{2L}{\kappa\,\mathrm{diam}(K)}$ of the $P$ DOFs will be EPWs.

In 3D, we propose a similar recipe with $L=\sqrt{\frac P8}$ and
$\zeta_p=\max\{1,\frac {2L}{\kappa\,\mathrm{diam}(K)}\sqrt{\xi_p}\}$.
In this case the uniform sampling of the orthogonal directions $\bp$ and $\be$ is more delicate and we refer to \cite{GaMoPa2026} for the details.

\section{Ultraweak variational formulation}\label{s:UWVF}

The EPW discrete space \eqref{eq:EPWspace} can be used within any of the Trefftz variational formulations in \cite[\S2]{TrefftzSurvey} (LS, TDG, UWVF, VTCR, WBM,\ldots).
Here we focus on the UWVF \cite{CED98,HMK02}, which can be recast as a special instance of the TDG \cite{HiMoPe11,TrefftzSurvey,GaMoPa2026}.

Consider the Helmholtz impedance BVP
\begin{align}\label{eq:BVP}\begin{aligned}
\Delta u +\kappa^2 u=0  \quad & \text{on }\Omega\subset\IR^n,\\
\dn u-\imath\kappa\sigma u=g \quad & \text{on }\partial\Omega,
\end{aligned}\end{align}
with $0<\sigma\in L^\infty(\deO)$ and $g\in L^2(\deO)$.
Let $\Omega$ be Lipschitz, bounded, partitioned in a mesh $\calT=\{K\}$, and let $\calF:=(\bigcup_{K\in\calT}\deK)\setminus\deO$ be the internal mesh skeleton.
Define the Trefftz space
\begin{align*}
T(\calT):=&\prod_{K\in\calT}T(K),\\
T(K):=&\big\{v\in H^1(K):\,\Delta v+\kappa^2 v=0,\\ &\qquad \dnK v\in L^2(\deK)\big\}.
\end{align*}
Extend $\sigma$ to a positive $L^\infty(\deO\cup\calF)$ and define two Robin traces:
$$
\gamma_\pm^K: T(K)\to L^2(\deK),\; \gamma_\pm^Kv:=\pm\dnK v-\imath\kappa\sigma v.
$$
For a finite dimensional subspace $V_h$ of $T(\calT)$, the classical UWVF \cite{CED98,HMK02} of BVP \eqref{eq:BVP} consists in finding $u_h\in V_h$ such that 
\begin{align}\label{eq:UWVF}
&\sum_{K\in\calT} \int_\deK \sigma^{-1}\,\gamma_-^K u_h\,\overline{\gamma_-^K v_h}\\
&-\sum_{K_1\in\calT}\sum_{K_2\in\calT} \int_{\deK_1\cap\deK_2}
  \sigma^{-1}\,\gamma_-^{K_1} u_h\, \overline{\gamma_+^{K_2} v_h}
 \nonumber\\
& =
  \sum_{K\in\calT}\int_{\deK\cap\deO} \sigma^{-1} \,g\,\overline{\gamma_+^K v_h}
\qquad \forall v_h\in V_h.
\nonumber
\end{align}
See \cite{GaMoPa2026} for the extension to piecewise-constant parameters and general mixed (Dirichlet, Neumann, impedance) boundary conditions.
The UWVF \eqref{eq:UWVF} is consistent with the BVP \eqref{eq:BVP}, well-posed, coercive and quasi-optimal in skeleton norms \cite{HiMoPe11,TrefftzSurvey,GaMoPa2026}.

The linear system of \eqref{eq:UWVF} can be written as $(\mD-\mC)\bu=\bb$, where $\mD$ and $\mC$ are the matrices corresponding to the two integrals at the left-hand side.
The matrix $\mD$ is Hermitian, positive-definite, and block-diagonal with each block corresponding to a mesh element, so the system is equivalent to $(\mI-\mD^{-1}\mC)\bu=\mD^{-1}\bb$ as in  \cite[eq.~(2.32)]{CED98}.
The inversion of $\mD$ is a local (elementwise) operation.

If the mesh facets are flat and EPW or PPW bases are used, then the entries of the matrices $\mD,\mC$ can be computed analytically.

\paragraph{Oversampled regularised UWVF}
We use the UWVF \eqref{eq:UWVF} with the sampled EPW discrete space from \eqref{eq:EPWspace}:
$$ 
V_h =\big\{v\in L^2(\Omega):\; v|_K\in \IT_P(K) \;\;\forall K\in\calT\big\}.
$$
The sampling of \S\ref{s:Recipe} can be done independently for each element $K$ or, if the element sizes are comparable, once for all elements.

However, following \cite{AdHu2019}, to ensure a stable solution even with basis
functions that are close to linearly dependent, we need to use
\textbf{oversampling} (10\% oversampling is used in \S\ref{s:Numerics}) and
\textbf{regularisation}.
The regularisation relies on the truncated SVD of the $\mD$ matrix, which can be done in parallel (elementwise) as in \cite{BaBeDiTo2021}.

We proceed as follows.
For each $K\in\calT$, let $\IT_{N\trial_K}(K)$ and $\IT_{N\test_K}(K)$ be two
EPW spaces in the form \eqref{eq:EPWspace}, with $N\test_K\ge N\trial_K$,
and define $N^\bullet:=\sum_{K\in\calT}N^\bullet_K$ for $\bullet\in\{$trial,test$\}$.
Denote $\mD_K\in \IC^{N\test_K\times N\trial_K}$ the tall matrix corresponding
to the diagonal block associated to the element \(K\) of the Galerkin matrix 
\(\mD\in \IC^{N\test\times N\trial}\) of the first integral in \eqref{eq:UWVF}.
Let $\mD_K=\mU_{K}\mSigma_{K}\mV_{K}^*$ be its SVD and
$\sigma_1\ge\sigma_2\ge\cdots\ge\sigma_{N_K\trial}\ge0$ be the singular values.
For a regularisation parameter $\epsilon\in(0,1)$ ($\epsilon=10^{-14}$ in
\S\ref{s:Numerics}), larger than the machine precision, let
$\mSigma_{K,\epsilon}^\dagger$ be the truncated pseudo-inverse of $\mSigma$:
the diagonal matrix with \(q\)-th diagonal entry
$\sigma_q^{-1}$ if $\sigma_q\ge\epsilon\sigma_1$ and $0$ otherwise.
We then define
$\mD_{\epsilon}^\dagger\in\IC^{N\trial\times N\test}$ the 
matrix with diagonal blocks $\mV_{K}\mSigma_{K,\epsilon}^\dagger\mU_{K}^*$.
Denoting by $\mC\in \IC^{N\test\times N\trial}$ the Galerkin matrix of the
second integral in \eqref{eq:UWVF},
the linear system $(\mI-\mD^{-1}\mC)\bu=\mD^{-1}\bb$ is approximated by the
$N\trial\times N\trial$ system
\begin{equation}\label{eq:LSEreg}
(\mI-\mD_\epsilon^\dagger\mC)\bu_\epsilon=\mD_\epsilon^\dagger\bb.
\end{equation}
The regularized linear system \eqref{eq:LSEreg} is sparse and can be solved by any direct or iterative method.

\section{Numerical experiments}\label{s:Numerics}

We finally illustrate the differences in the approximation capabilities of PPWs (taking \(L=0\) in the sampling recipe in \S\ref{s:Recipe}) and EPWs. 

\paragraph{Point source}

The domain is the square \(\Omega:=(0,1) \times (-\frac12,\frac12)\), discretized into an unstructured quasi-uniform mesh with \(41\) triangles.
We consider the numerical approximation of the fundamental solution
\(u(\mathbf{x}) = \frac\imath4 H^{(1)}_{0}(\kappa |\mathbf{x} - \mathbf{s}|)\)
using the UWVF~\eqref{eq:UWVF} with a suitable datum \(g\).
The wavenumber is \(\kappa=128\) and the singularity \(\mathbf{s}=(-\frac{\pi}{5\kappa},0)\)
is located a tenth of a wavelength away from the left boundary.

\begin{figure}[htbp]
    \centering
    \includegraphics[width=0.97\linewidth]{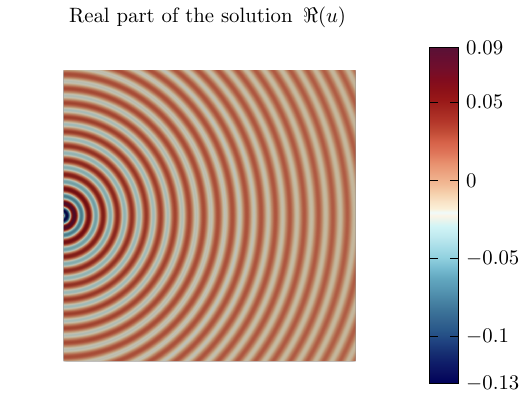}
    \includegraphics[width=0.97\linewidth]{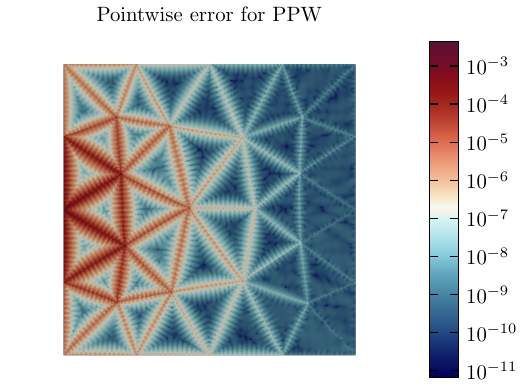}
    \includegraphics[width=0.97\linewidth]{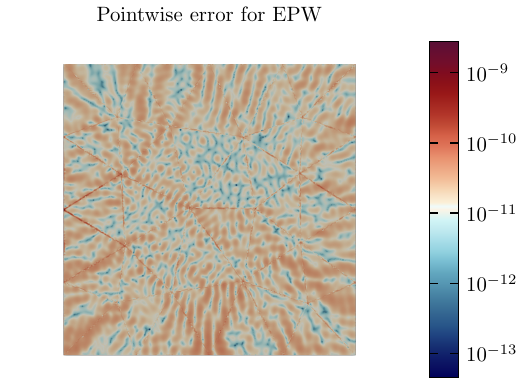}
    \includegraphics[width=\linewidth]{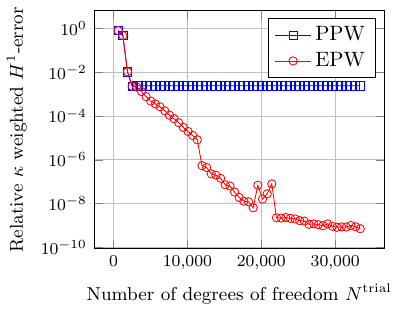}
    \caption{Fundamental solution:
        \(\Re({u})\) (top),
        pointwise error (middle) for PPWs and EPWs,
        convergence history (bottom).
    }
    \label{fig:hankel}
\end{figure}

Figure~\ref{fig:hankel} reports the real part of the solution $u$, the convergence of the \(\kappa\)-weighted relative \(H^{1}(\Omega)\)-error with respect to the number of degrees of freedom, and the pointwise error in the domain (for the largest approximation space, $P=815$) when using either PPWs or EPWs.
Note the different scales in the error colorbars.

We observe a first regime of fast convergence, obtained for both types of waves
since the two approximation spaces are actually of similar nature (propagative)
when the number of degrees of freedom is moderate.
Then the PPW instability makes convergence stall and prevents from obtaining higher accuracy, whereas EPWs incur a much smaller error.
The ratio between the EPW and the PPW errors reaches values below $10^{-6}$.
The convergence is not monotone since the approximation spaces are not nested.

\paragraph{Influence of the wavenumber}

In the same configuration, we now study the influence of the wavenumber \(\kappa\) when the size of the approximation space is proportional to \(\kappa\) itself.
More precisely, the number of PWs per cell is fixed to \(P=4\kappa\).
As before, the singularity is placed at \(\mathbf{s}=(-\frac{\pi}{5\kappa},0)\), lying at a distance from the boundary proportional to the wavelength.

\begin{figure}[htbp]
 \centering
 \includegraphics[width=\linewidth]{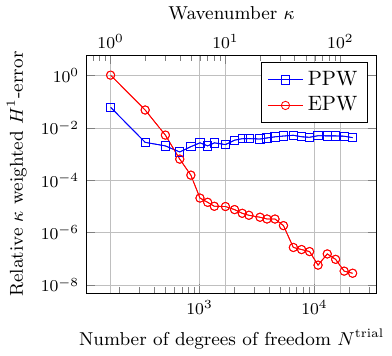}
 \caption{Fundamental solution: convergence result for increasing 
  \(\kappa\) using an approximation space of size increasing linearly in
  \(\kappa\).}
  \label{fig:vs_k}
\end{figure}

Figure~\ref{fig:vs_k} reports the convergence of the \(\kappa\)-weighted relative \(H^{1}(\Omega)\)-error with respect to the number of degrees of freedom and the wave\-number.
The error achieved using EPWs improves with the wavenumber as the size of the approximation space increases linearly with~\(\kappa\).

\paragraph{Scattering problem}

We consider the scattering of an incident plane wave
\(u_{i}(\mathbf{x}) := e^{\imath \kappa \mathbf{d} \cdot \mathbf{x}}\)
for \(\mathbf{d} = (\frac12,-\frac{\sqrt3}2)\) and wavenumber \(\kappa=16\) by a sound-soft (i.e.\ Dirichlet) non-convex scatterer of diameter \(2\).
The UWVF formulation in the total field features a non-homogeneous Robin boundary condition on a truncation boundary of diameter \(4\) surrounding the obstacle.
The mesh contains 64 triangles.
Since the solution to this problem is not known analytically, we use as reference solution the EPW approximation with twice as many waves as the largest approximation space.
This test case is more challenging because of the presence of the cavity and the corner singularities at each vertex of the polygonal scatterer.

Figure~\ref{fig:scattering} reports the real part of the solution, the convergence history of the \(\kappa\)-weighted relative \(H^{1}(\Omega)\)-error with respect to the number of degrees of freedom and the pointwise error in the domain (for the largest approximation space, $P=2080$) when using either PPWs or EPWs.

Again we observe that the convergence of the UWVF using PPWs stagnates whereas EPWs yield much better accuracy.

\begin{figure}[htbp]
    \centering
    \includegraphics[width=0.9\linewidth]{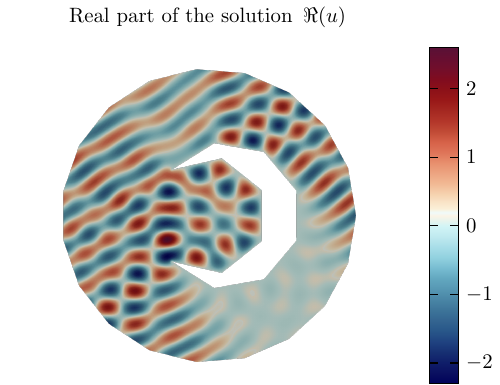}
    \includegraphics[width=0.9\linewidth]{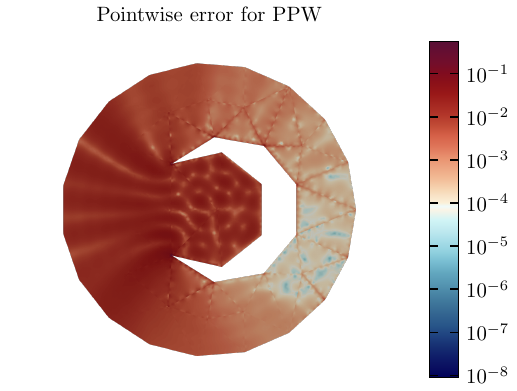}
    \includegraphics[width=0.9\linewidth]{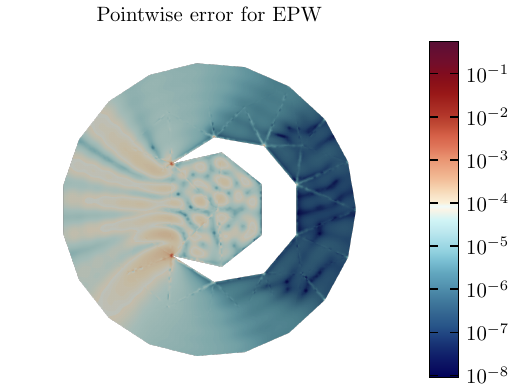}
    \includegraphics[width=\linewidth]{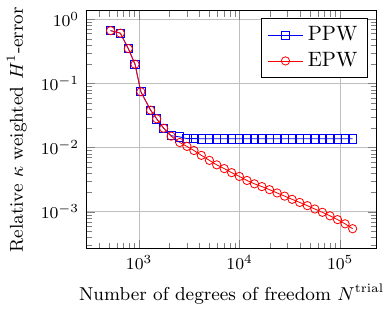}
    \caption{Scattering problem:
           \(\Re({u})\) (top),
           pointwise error (middle) for PPWs and EPWs,
           convergence history (bottom).
    }
    \label{fig:scattering}
\end{figure}

\end{wavespaper}

\begin{thebibliography}{99}
\bibitem{AdHu2019}
B. Adcock, and D. Huybrechs, Frames and numerical approximation. \textit{SIAM Rev.} \textbf{61} (2019), 443--473.

\bibitem{BaBeDiTo2021}
H. Barucq, A. Bendali, J. Diaz, and S. Tordeux,
Local strategies for improving the conditioning of the plane-wave ultra-weak variational formulation. 
\textit{J. Comput. Phys.} \textbf{441} (2021),  110449, 18.

\bibitem{CED98}
O. Cessenat, and B. Després,
Application of an ultra weak variational formulation of elliptic PDEs to the two-dimensional Helmholtz equation. 
\textit{SINUM} \textbf{35.1} (1998), 255--299.

\bibitem{CoNi2025}
J. Coyle, and N. Nigam, 
The whys and hows of conditioning of DG plane wave Trefftz methods: a single element.\\
\textit{arXiv:2509.14500} (2025).

\bibitem{Galante2025}
N. Galante, A. Moiola, and E. Parolin,
Stable approximation of Helmholtz solutions in the 3D ball using evanescent plane waves. 
\textit{SMAI J. Comput. Math.}
\textbf{11} (2025), 435--472.

\bibitem{GaMoPa2026}
N. Galante, A. Moiola, and E. Parolin,
Trefftz Discontinuous Galerkin method using evanescent plane waves.
In preparation (2026).

\bibitem{HiMoPe11}
R. Hiptmair, A. Moiola, and I. Perugia,
Plane wave discontinuous Galerkin methods for the 2D Helmholtz equation: analysis of the p-version. 
\textit{SINUM} \textbf{49} (2011), 264--284.

\bibitem{TrefftzSurvey}
R. Hiptmair, A. Moiola, and I. Perugia,
A survey of Trefftz methods for the Helmholtz equation,
\textit{Lect. Notes Comput. Sci. Eng.} (2016), 237--278.

\bibitem{HMK02}
T. Huttunen, P. Monk, and J. P. Kaipio,
Computational aspects of the ultra-weak variational formulation. 
\textit{J. Comput. Phys.} \textbf{182.1} (2002), 27--46.

\bibitem{MoHiPe11}
A. Moiola, R. Hiptmair, and I. Perugia,
Plane wave approximation of homogeneous Helmholtz solutions. 
\textit{Z.Angew. Math. Phys.} \textbf{62} (2011), 809--837.

\bibitem{Parolin2023}
E. Parolin, D. Huybrechs, and A. Moiola,
Stable approximation of Helmholtz solutions in the disk by evanescent plane waves.
\textit{M2AN} \textbf{57(6)} (2023), 3499--3536.

\bibitem{Trefftz1926}
E. Trefftz, Ein Gegenstuck zum Ritzschen Verfahren, 
\textit{Proceedings of the 2nd International Congress for Applied Mechanics}, Zurich (1926), 131--137.

\end{thebibliography}
\end{document}